\documentclass[12pt]{article}

\usepackage{typearea}
\usepackage{a4wide}
\usepackage{amsmath}
\usepackage{amsfonts}
\usepackage[frenchb, english]{babel}
\usepackage[T1]{fontenc} 
\usepackage[dvips
]{graphicx}
\usepackage{mathrsfs}

\usepackage{aeguill}
\usepackage{eufrak}

\usepackage{latexsym,amsmath,amssymb,amsfonts,amsthm,verbatim}
\usepackage[utf8]{inputenc}
\usepackage{color}
\usepackage{eufrak}
\usepackage{cases}

\let \ssection=\section
\renewcommand{\section}{\setcounter{equation}{0}\ssection}

\newtheorem{Th}{Theorem}[section]
\newtheorem{theorem}[Th]{Theorem}

\newtheorem{proposition}[Th]{Proposition}
\newtheorem{lemma}[Th]{Lemma}
	\theoremstyle{definition}

	\theoremstyle{remark}
\newtheorem{remark}[Th]{Remark}
	\theoremstyle{example}

\newcommand{\suppX}{supp(X_t)}

\newcommand{\Dtht}{D^2_{\theta,t}}

\newcommand{\Dthdeux}{D_{\theta}}
\newcommand{\Dtun}{D_{t}}
\newcommand{\pbg}{\left(\beta\circ g^{-1}\right)}

\newcommand{\fxx}{f_{xx}\lp X_s,Y_s\rp}
\newcommand{\fyy}{f_{yy}\lp X_s,Y_s\rp}
\newcommand{\fyx}{f_{yx}\lp X_s,Y_s\rp}

\newcommand{\fyyr}{f_{yy}\lp X_r,Y_r\rp}
\newcommand{\fyxr}{f_{yx}\lp X_r,Y_r\rp}

\newcommand{\Dt}{D_{\theta}}
\newcommand{\Dtt}{D^2_{\theta,t}}
\newcommand{\xys}{\lp X_s,Y_s\rp}
\newcommand{\xyr}{\lp X_r,Y_r\rp}

\newcommand{\E}{\mathbf{E}}

\newcommand{\F}{\mathcal F}

\newcommand{\lp}{\left (}
\newcommand{\rp}{\right )}

\newcommand{\undeux}{\frac{1}{2}}

\newcommand{\intt}{\int_t^T}
\newcommand{\into}{\int_0^t}

\newcommand{\intts}{\int_t^s}
\newcommand{\inttht}{\int_\theta^t}

\newcommand{\dxt}{D_{\theta}X_t}
\newcommand{\dut}{D_{\theta}U_t}
\newcommand{\dxs}{D_{\theta}X_s}
\newcommand{\dus}{D_{\theta}U_s}
\newcommand{\dyt}{D_{\theta}Y_t}
\newcommand{\dyT}{D_{\theta}Y_T}
\newcommand{\dys}{D_{\theta}Y_s}
\newcommand{\dzt}{D_{\theta}Z_t}
\newcommand{\dzs}{D_{\theta}Z_s}

\newcommand{\xyzt}{\lp X_t,Y_t,Z_t\rp}
\newcommand{\xyzs}{\lp X_s,Y_s,Z_s\rp}
\newcommand{\xyzr}{\lp X_r,Y_r,Z_r\rp}
\newcommand{\fxs}{f_x\xyzs}

\newcommand{\fxt}{f_x\xyzt}
\newcommand{\fys}{f_y\xyzs}
\newcommand{\fyr}{f_y\xyzr}
\newcommand{\fyt}{f_y\xyzt}
\newcommand{\fzs}{f_z\xyzs}

\newcommand{\fzt}{f_z\xyzt}

\newcommand{\Ept}[1]{\widetilde{\mathbf{E}}\lp\left. #1 \right| \F_t\rp}
\newcommand{\Et}[1]{\E\lp\left. #1 \right| \F_t\rp}

\begin{document}


\renewcommand{\thefootnote}{\fnsymbol{footnote}}

\renewcommand{\thefootnote}{\fnsymbol{footnote}}
\begin{center}
{\large{ \textsc{Density estimates for solutions to one dimensional Backward SDE's}}}\\~\\
Omar Aboura\footnote{
SAMM, EA 4543, Universit\'e Paris 1 Panth\'eon Sorbonne, 90 Rue de Tolbiac, 75634 Paris Cedex France.
Email: {\tt omar.aboura@malix.univ-paris1.fr}
}
and 
Solesne Bourguin\footnote{
SAMM, EA 4543, Universit\'e Paris 1 Panth\'eon Sorbonne, 90 Rue de Tolbiac, 75634 Paris Cedex France.
}\footnote{
Faculté des Sciences, de la Technologie et de la Communication; UR en Mathématiques. 6, rue Richard Coudenhove-Kalergi, L-1359 Luxembourg.
Email: {\tt solesne.bourguin@uni.lu}
}
\end{center}
{\small \noindent {\bf Abstract:} 
In this paper, we derive sufficient conditions for each component of the solution to a general backward stochastic differential equation to have a density for which upper and lower Gaussian estimates can be obtained.
\\


\noindent
{\bf 2010 AMS Classification Numbers:} 
60H10,
60H07.
\\


\section{Introduction}
\noindent In \cite{NV}, I. Nourdin and F.G. Viens have introduced sufficient conditions to prove the existence of a density for a Malliavin differentiable random variable and to provide upper and lower Gaussian estimates for this density.

\noindent This result has led to several research papers, such as those by D. Nualart and L.~Quer-Sardanyons (\cite{NQ1}, \cite{NQ2}), in which these authors applied Nourdin and Viens result to solutions of quasi-linear stochastic partial differential equations and to a class of stochastic equations with additive noise.

\noindent In this paper, we use Nourdin and Viens's approach to prove that, under proper conditions on the coefficients, each component of the solution $(X_t,Y_t,Z_t)$ to a backward stochastic differential 
equation
\begin{numcases}
~X_t = x_0 + \into b(X_s) ds + \into \sigma(X_s) dW_s \label{diffusion}\\
Y_t = \xi + \intt f \xyzs ds - \intt Z_s dW_s       \label{backSDE}
\end{numcases}
has a density for which upper and lower Gaussian bounds can be derived. This implies to study the relation between the coefficients of the diffusion equation (\ref{diffusion}) and the coefficients of the backward SDE (\ref{backSDE}).

\noindent Our paper is organized as follows: in a first part, we study the component $Y_{t}$ of the solution and in a second and last part, we focus on the component $Z_{t}$ (for which, up to our knowledge, no density existence results exist). We will not develop a specific study of the first component $X_{t}$ of the solution to the BSDE due to the fact that the question of the existence of a density for the solution to an SDE of the type (\ref{diffusion}) and the properties of this density has been intensively studied and we refer the reader to \cite{N} for an extensive survey of the existing litterature and results on this topic. 

\noindent Equations of the type (\ref{backSDE}) were introduced in \cite{PP90} and are closely related with viscosity solution to PDEs. These equations have been intensively studied and have many applications in control theory and financial methods among others.

\noindent The existence of the density for the random variable $Y_t$ at a fixed time $t\in(0,T)$, as well as upper bounds for its tail behavior, have been proven by F. Antonelli and A. Kohatsu-Higa \cite{AKH}, using the Bouleau-Hirsch Theorem. We retrieve Antonelli and Kohatsu-Higa's existence result for the density of $Y_t$, and we also derive Gaussian estimates for it. In order to provide (additionally to the existence result itself) estimates for the density of $Y_t$, we need to strengthen the hypotheses of Antonelli and Kohatsu-Higa. 

\noindent We also address the question of the existence of a density for the random variable $Z_t$ as well as the possibility of deriving Gaussian estimates for it. This question has not been solved in \cite{AKH}. We need the same hypotheses as in the case of $Y_t$, as well as additional ones, since $Z_t$ can be expressed as a function of the Malliavin derivative of $Y_t$.

\noindent In order to be self contained, we first give an overview of some elements of Malliavin calculus in Section 2, and the corresponding notations. Section 3 is dedicated to the component $Y_{t}$ of the solution to the BSDE. Section 4 deals with the question of the existence of a density for $Z_{t}$ and is organized in two subsections, dealing respectively with the question of the existence of a density and conditions for this density to be bounded by Gaussian upper and lower estimates.

\section{Framework, main tools and notations}

\subsection{Elements of Malliavin calculus}

\noindent Consider the real separable Hilbert space $L^{2}(\left[0,T\right])$ and let $(W (\varphi), \varphi\in L^{2}(\left[0,T\right]))$ be an isonormal Gaussian process 
on a probability space $(\Omega, {\EuFrak{A}}, P)$, 
that is, a centered Gaussian family of random variables such that $\mathbf{E}\left( W(\varphi) W(\psi) \right)  = \langle\varphi, \psi\rangle_{L^{2}(\left[0,T\right])}$. For any integer $n \geq 1$, denote by $I_{n}$ the multiple stochastic integral with respect to $W$ (see \cite{N} for an extensive survey on Malliavin calculus). The map $I_{n}$ is actually an isometry between the Hilbert space $L^{2}(\left[0,T\right]^{n})$ equipped with the scaled norm $\frac{1}{\sqrt{n!}}\Vert\cdot\Vert_{L^{2}(\left[0,T\right]^{n})}$ and the Wiener chaos of order $n$, which is defined as the closed linear span of the random variables $H_{n}(W(\varphi))$ where $\varphi\in L^{2}(\left[0,T\right]), \Vert\varphi\Vert_{L^{2}(\left[0,T\right])}=1$ and $H_{n}$ is the Hermite polynomial of degree $n\geq 1$, that is defined by
\begin{equation*}
H_{n}(x)=\frac{(-1)^{n}}{n!} \exp \left( \frac{x^{2}}{2} \right)
\frac{d^{n}}{dx^{n}}\left( \exp \left( -\frac{x^{2}}{2}\right)
\right), \hskip0.5cm x\in \mathbb{R}.
\end{equation*}
The isometry of multiple integrals can be written as follows: for positive integers $m,n$,
\begin{eqnarray*}
\mathbf{E}\left(I_{n}(f) I_{m}(g) \right) &=& n! \langle f,g\rangle _{L^{2}(\left[0,T\right]^{n})}\quad \mbox{if } m=n,\nonumber \\
\mathbf{E}\left(I_{n}(f) I_{m}(g) \right) &= & 0\quad \mbox{if } m\not=n.
\end{eqnarray*}
It also holds that
\begin{equation*}
I_{n}(f) = I_{n}\big( \tilde{f}\big)
\end{equation*}
where $\tilde{f} $ denotes the symmetrization of $f$ defined by $$\tilde{f}%
(x_{1}, \ldots , x_{n}) =\frac{1}{n!} \sum_{\sigma \in {\mathcal S}_{n}}
f(x_{\sigma (1) }, \ldots , x_{\sigma (n) } ). $$
We recall that any square integrable random variable $F$ which is measurable with respect to the $\sigma$-algebra generated by $W$ can be expanded into an orthogonal sum of multiple stochastic integrals
\begin{equation}
\label{sum1} F=\sum_{n\geq0}I_{n}(f_{n})
\end{equation}
where $f_{n}\in L^{2}(\left[0,T\right]^{n})$ are (uniquely determined)
symmetric functions and $I_{0}(f_{0})=\mathbf{E}\left[  F\right]$.
\\\\
\noindent Let $L$ be the Ornstein-Uhlenbeck operator defined by $LF=-\sum_{n\geq 0} nI_{n}(f_{n})$ if $F$ is given by (\ref{sum1}) and satisfies $\sum_{n \geq 1}n^{2}n!\Vert f_{n} \Vert^{2} < \infty$. For $p>1$ and $\alpha \in \mathbb{R}$ we introduce the Sobolev-Watanabe space $\mathbb{D}^{\alpha ,p }$  as the closure of
the set of random variables of the form (\ref{vasimple}) (see (1.28) in \cite{N}) with respect to the norm defined by
\begin{equation*}
\Vert F\Vert _{\alpha , p} =\Vert (I -L) ^{\frac{\alpha }{2}}F \Vert_{L^{p} (\Omega )},
\end{equation*}
where $I$ represents the identity. We denote by $D$  the Malliavin  derivative operator that acts on smooth random variables of the form 
\begin{equation}
\label{vasimple} F=g(W(\varphi _{1}), \ldots , W(\varphi_{n})),
\end{equation}
where $g$ is a smooth function with compact support and $\varphi_{i} \in  L^{2}(\left[0,T\right])$, as follows:
\begin{equation*}
DF=\sum_{i=1}^{n}\frac{\partial g}{\partial x_{i}}(W(\varphi _{1}), \ldots , W(\varphi_{n}))\varphi_{i}.
\end{equation*}
The operator $D$ is continuous from $\mathbb{D} ^{\alpha , p} $ into $\mathbb{D} ^{\alpha -1, p} \left( L^{2}(\left[0,T\right])\right).$ The adjoint of $D$ is denoted by $\delta $ and is called the divergence (or Skorohod) integral. It is a continuous operator from $\mathbb{D}^{\alpha, p } \left( {L^2([0,T])}\right)$ into $\mathbb{D} ^{\alpha +1,p}$. More generally, we can introduce iterated weak derivatives of order $k$. If $F$ is a smooth random variables and $k$ is a positive integer, we set $$D_{t_{1},...,t_{k}}^{k}F = D_{t_{1}}D_{t_{2}}...D_{t_{k}}F.$$ 
We have the following duality relationship between $D$ and $\delta$ for $F \in \mathbb{D}^{1,2}$ and $u \in dom$ $\delta$
\begin{equation*}
\mathbf{E} (F\delta (u))= \mathbf{E}\langle DF, u\rangle _{L^2([0,T])}.
\end{equation*}
For adapted integrands, the divergence integral coincides with the classical It\^o integral. We will use the notation
\begin{equation*}
\delta (u) =\int_{0}^{T} u_{s} dW_{s}.
\end{equation*}
Note that the following integration by parts relation between $D$ and $\delta$ holds
\begin{equation*}
D_{t}(\delta (u)) = u_{t} + \int_0^{T}D_{t}u_{s}dW_{s},
\end{equation*}
where $u \in \mathbb{D}^{1,2}(L^{2}([0,T]))$ is such that $\delta (u) \in \mathbb{D}^{1,2}$.

\subsection{Density existence and Gaussian estimates}

\noindent A classical density existence result is the celebrated Bouleau-Hirsch theorem (see \cite{N} for an extensive survey on this result). This result provides conditions in term of Malliavin derivatives for a random variable to have a density.
\begin{theorem}[Bouleau--Hirsch]
\label{bouleauhirsch}
\noindent Let $F$ be a random variable of the space $\mathbb{D}^{1,2}$ and suppose that $\Vert DF \Vert _{L^{2}(\left[0,T\right])} > 0$ a.s. Then the law of $F$ is absolutely continuous with respect to the Lebesgue measure on $\mathbb{R}$.
\end{theorem}

\noindent In \cite{NV}, Corollary 3.5, Nourdin and Viens have given the following sufficient condition for a weakly differentiable random variable to have a density with lower and upper Gaussian estimates.
\begin{proposition}
\label{existenceBornes}
Let $F$ be in $\mathbb{D}^{1,2}$ and let the function $g$ be defined for all $x \in \mathbb{R}$ by
\begin{align}
\label{g2z}
g(x) = \mathbf{E}\left( \langle DF,-DL^{-1}F \rangle_{L^{2}(\left[0,T\right])}\Big|F - \mathbf{E}(F) = x\right).
\end{align}
If there exist positive constants
$\gamma_{\min},\gamma_{\max}$ such that, for all $x \in \mathbb{R}$, almost surely
\[
0<\gamma_{\min}^{2}\leqslant g(x)  \leqslant\gamma_{\max}^{2}%
\]
then $F$ has a density $\rho$ satisfying, for almost all $z\in\mathbb{R}$%
\[
\frac{\mathbf E|F- \mathbf{E}(F)|}{2\,\gamma_{\max}^{2}}\,\mathrm{exp}\left(  -\frac{(z- \mathbf{E}(F))^{2}}%
{2\gamma_{\min}^{2}}\right)  \leqslant\rho(z)\leqslant\frac{\mathbf E|F- \mathbf{E}(F)|}%
{2\,\gamma_{\min}^{2}}\,\mathrm{exp}\left(  -\frac{(z- \mathbf{E}(F))^{2}}{2\gamma_{\max}^{2}%
}\right)  .
\]
\end{proposition}
\noindent Furthermore, Nourdin and Viens have also provided the following useful result, which gives some rather explicit description of $g(x)$.
Recall that $W=\lp W(\phi),\phi\in L^2\lp[0,T]\rp\rp$.
\begin{proposition}
\label{PropPourCalcg}
Let $F$ be in $\mathbb{D}^{1,2}$ and write $DF=\Phi_{F}(W)$ with a measurable function
$\Phi_{F}:\mathbb{R}^{L^{2}(\left[0,T\right])}\rightarrow L^{2}(\left[0,T\right])$. Then, if $g(x)$ is defined by (\ref{g2z}), we have
\[
g(x)  =\int_{0}^{\infty}e^{-u}\,\mathbf{E}\left(\mathbf{E'}\big(\langle\Phi
_{F}(W), \widetilde{\Phi
_{F}^{u}}(W) \rangle_{L^{2}(\left[0,T\right])}\big)|F - \mathbf{E}(F) = x\right)du,
\]
where $\widetilde{\Phi_{F}^{u}}(W) = \Phi_{F}(e^{-u}W+\sqrt{1-e^{-2u}}W^{\prime})$, $W^{\prime}$ stands for an independent copy of $W$, and is such that $W$
and $W^{\prime}$ are defined on the product probability space $(\Omega
\times\Omega^{\prime},\F\otimes\F^{\prime},\mathbb{P}\times
\mathbb{P}^{\prime})$ and $\mathbf{E'}$ denotes the mathematical expectation with
respect to $\mathbb{P}^{\prime}$.
\end{proposition}

\subsection{Notations}

We denote by $\mathcal{C}_{b}^{n}(\mathbb{R}^{p})$ the space of $n$--times differentiable functions on $\mathbb{R}^{p}$ with bounded partial derivatives up to order $n$.
\\~\\
Let $f$ be a three times differentiable function of three variables $x$, $y$ and $z$. We will use the following notations : $\frac{\partial f}{\partial x} = f_{x}$, $\frac{\partial f}{\partial y} = f_{y}$, $\frac{\partial f}{\partial z} = f_{z}$ $\frac{\partial^{2} f}{\partial x^{2}} = f_{xx}$, $\frac{\partial^{2} f}{\partial y^{2}}(x,y) = f_{yy}$, $\frac{\partial^{2} f}{\partial x \partial y} = f_{xy}$, 
$\frac{\partial^{2} f}{\partial y \partial x} = f_{yx}$. \\~\\We will also use the following notation for the Lie bracket : $[h,g]=hg'-gh'$, where $h,g: \mathbb{R} \rightarrow \mathbb{R}$.
~\\~\\
In the whole paper, $c$ and $C$ will denote constants that may vary from line to line.

\section{Density of $Y_t$ : existence and Gaussian estimates}

\noindent The following backward stochastic differential equation was introduced in Pardoux and Peng \cite{PP90} (see also \cite{PP}) and density properties of its solutions were investigated in \cite{AKH}:
\begin{numcases}
~X_t = x_0 + \into b(X_s) ds + \into \sigma(X_s) dW_s \label{eqdiff}\\
Y_t = \xi + \intt f \xyzs ds - \intt Z_s dW_s       \label{eqbackSDE}
\end{numcases}
In this section, we give conditions for the random variable $Y_t$ to have a density which can be bounded from above and below by Gaussian ones.

\label{subsectionHypY}
\subsection{Hypotheses}\label{hyp-Y}
We consider $b$, $\sigma$ and $f$ to be appropriately smooth functions to ensure the existence and uniqueness of solutions to equations (\ref{eqdiff}) and (\ref{eqbackSDE}). We also impose additional conditions needed to state our main result:
\begin{numcases}
~\mbox{\textbf{H1 : }} \xi \in L^{2}(\Omega, \mathcal{F}_{T})\cap \mathbb{D}^{1,2} \quad \mbox{and} \quad \forall \mbox{\ } \theta \leq T,\mbox{\ } 0 < c \leq D_{\theta}\xi \leq C \quad a.s.  \nonumber \\
\mbox{\textbf{H2 : }} f \in \mathcal{C}_{b}^{1}(\mathbb{R}^{3}) \quad \mbox{and} \quad 0 \leq f_{x} \leq C. \nonumber \\
\mbox{\textbf{H3 : }} b \in \mathcal{C}_{b}^{1}(\mathbb{R}), \quad \sigma \in \mathcal{C}_{b}^{2}(\mathbb{R}), \quad 0 \leq \sigma \leq C  \quad \mbox{and} \quad \vert[b,\sigma]\vert \leq M\sigma. \nonumber
\end{numcases}
where $[b,\sigma]$ denotes the Lie bracket between $b$ and $\sigma$. 
\begin{remark}
Note that it is natural to have a condition on the Lie bracket between $b$ and $\sigma$ as this quantity is the one that appearing in the classical density and smoothness results (e.g. Hörmander's brackets condition). 
\end{remark}
\begin{remark}
The hypotheses on the positivity of $\sigma$ and $f_{x}$ are made in order to make the proofs as readable as possible. In fact, one only needs $\sigma$ and $f_{x}$ to have the same sign to draw to same conclusions.
\end{remark}

\label{mainressec}
\subsection{Main result}\label{mainres}
\noindent Under the above assumptions, $Y_t$ has a density for which the following Gaussian estimates can be derived.
\begin{theorem}
\label{theoY}
Under the above hypotheses, for $t\in(0,T)$ the random variable $Y_t$ defined in (\ref{eqbackSDE}) has a density $\rho_{Y_t}$. Furthermore, there exist strictly positive constants $c$ and $C$ such that, for almost all $y\in\mathbb{R}$ and all $t \in \left[0,T \right]$, $\rho_{Y_t}$ satisfies the following:
\begin{align*}
\frac{\mathbf{E}\vert Y_t -\mathbf{E}(Y_t)\vert}{2ct}\mathrm{exp}\left(-\frac{(y-\mathbf{E}(Y_t))^{2}}{2Ct}\right)  
\leq \rho_{Y_t}(y) \leq \frac{\mathbf{E}\vert Y_t -\mathbf{E}(Y_t)\vert}{2Ct}\mathrm{exp}\left(-\frac{(y-\mathbf{E}(Y_t))^{2}}{2ct}\right).
\end{align*}
\end{theorem}

\noindent The rest of this section is devoted to the proof of Theorem \ref{theoY} which is divided in three steps. In the first step, we will prove that the Malliavin derivative of $X_{t}$ is bounded and non-negative. This property will be necessary in Step 2 (as $DX_{t}$ appears in the Malliavin derivative of $Y_{t}$).
\\~\\
In the second step, we will derive upper and lower bounds for the Malliavin derivative of $Y_{t}$. Indeed, our purpose is to use the Nourdin--Viens formula (Proposition \ref{existenceBornes}) in which one needs to bound a function of the Malliavin derivative of random variable for which density results are investigated. We will make use of the properties of $DX_{t}$ prooved in Step 1.
\\~\\
In the third step, we will use the bounds obtained on $DY_{t}$ and the Nourdin--Viens formula to conclude the proof.
\\~\\
\noindent {\bf Proof of Theorem \ref{theoY}: }
\\~\\
\textbf{\textit{Step 1: Boundedness and positivity of $DX_{t}$}}
\\~\\
Consider equation (\ref{eqdiff}). Using a Lamperti transform (see \cite{Lam64} or \cite{KS91} pp. 294--295 exercise 2.20), we compute the Malliavin derivative of $X_{t}$. The Lamperti transform of $X_{t}$, hereafter denoted by $U_{t}$, is given by
\begin{align*}
U_t = g(x_0) + \int_{0}^{t}\beta \circ g^{-1}(U_s)ds + W_t,
\end{align*}
where $$g(x) = \int_{0}^{x}\frac{du}{\sigma(u)} \quad \mbox{and} \quad \beta(x) =\frac{b}{\sigma}(x)- \frac{\sigma'(x)}{2}. $$ Computing the Malliavin derivative of $U_{t}$ yields, for $\theta \in \left[0,t \right]$, 
\begin{align}
\label{derU}
\dut = 1 + \inttht (\beta \circ g^{-1})'(U_s)\dus ds 
= \mbox{exp}\left[ \inttht (\beta \circ g^{-1})'(U_s) ds\right].
\end{align}
Deriving the identity $g \circ g^{-1}(x) = x$ on $g(\suppX)$ yields $(g^{-1})'(x) = \sigma \circ g^{-1}(x)$. Using this fact we get $(\beta \circ g^{-1})'(x) = \beta' \circ g^{-1}(x)(g^{-1})'(x) = (\beta' \sigma) \circ g^{-1}(x)$. In addition, it is easy to check that on $g(\suppX)$,
\begin{align}
\label{defBetaPrimeSigma}
(\beta' \sigma)(x) = \frac{[\sigma,b](x)}{\sigma(x)} - \frac{(\sigma \sigma'')(x)}{2}.
\end{align}
Gathering those results and using hypothesis $(\mbox{\textbf{H3}})$ of Subsection \ref{hyp-Y} immediately yields on $g(\suppX)$
\begin{align*}
-C \leq (\beta \circ G^{-1})' \leq C,
\end{align*}
where $C$ is a positive constant. Using (\ref{derU}), we deduce, $\mathbb{P}$-a.s,
\begin{align}
\label{BornesDU2}
0 < c \leq D_{\theta}U_t \leq C.
\end{align}
Furthermore, as $X_t = g^{-1}(U_t)$, it holds that, for $0 < \theta < t \leq T$,
\begin{align}
\label{relationEntreDXetDU2}
\dxt = (g^{-1})'(U_t)\dut = \sigma \circ g^{-1}(U_t)\dut.
\end{align}
Combining (\ref{BornesDU2}) and (\ref{relationEntreDXetDU2}) with the fact that $\sigma$ is bounded and non-negative yields, $\mathbb{P}$-a.s,
\begin{align}
\label{bornesDerXPreuvePropY}
0 \leq \dxt \leq C.
\end{align}
\\~\\
\textbf{\textit{Step 2: Computation of bounds on $DY_{t}$}}
\\~\\
We at first represent $\dyt$ by means of an equivalent probability; 
this is similar to \cite{AKH} and the proof is included for the sake of completeness. 
It is well known (see for example Theorem 2.2 in \cite{AKH}) that, for every $t\in(0,T]$, $Y_t \in \mathbb{D}^{1,2}$ and $Z \in L^{2}\left(0,T;\mathbb{D}^{1,2} \right)$. Furthermore, since $\theta < t$, we have 
\begin{align}
\label{numeroDyt}
\dyt =& D_{\theta}\xi -\intt \dzs dW_s \nonumber \\
&+\intt \left[ \fxs\dxs+\fys\dys+\fzs\dzs \right] ds.
\end{align}
The product $e^{\int_0^t \fys ds}\dyt$ yields a more suitable representation of $\dyt$; indeed, for $t \in (0,T]$,
and $0\leq\theta<t$
\begin{align*}
d\left[e^{\int_0^t \fys ds}\dyt\right] =& \Bigl[ \dyt e^{\int_0^t \fys ds} \fyt \Bigr. \\
& \Bigl.-  e^{\int_0^t \fys ds}\left(\fxt\dxt + \fyt\dyt\right.\Bigr. \\
& \Bigl.\left.+ \fzt\dzt\right)\Bigr] dt + e^{\int_0^t \fys ds}\dzt dW_t.
\end{align*}
Integrating from $t$ to $T$ yields, for $\theta < t$,
\begin{align*}
 e^{\int_0^T \fys ds}\dyT - &e^{\int_0^t \fys ds}\dyt = -\intt e^{\int_0^s \fyr dr}\left[\fxs\dxs \right. \\
& \left. + \fzs \dzs \right]ds  + \intt e^{\int_0^s \fyr dr}\dzs dW_s.
\end{align*}
Note that $\dyT = D_{\theta}\xi$; therefore, for $t \in (0,T]$,
\begin{align*}
\dyt =& e^{\intt \fys ds}D_{\theta}\xi + \intt e^{\intts \fyr dr}\left[\fxs\dxs \right. \\
& \left. + \fzs \dzs \right]ds  - \intt e^{\intts \fyr dr}\dzs dW_s.
\end{align*}
Let $\widetilde{W}_t = W_t - \int_{0}^{t}\fzs ds$. Because $f_z \in \mathcal{C}_{b}^{0}(\mathbb{R})$, Novikov's condition is verified and $\widetilde{W}$ is a Brownian motion under some equivalent probability $\widetilde{\mathbb{P}}$. Girsanov's theorem yields 
\begin{align*}
\dyt =& e^{\intt \fys ds}D_{\theta}\xi + \intt e^{\intts \fyr dr}\fxs\dxs ds \\
&  - \intt e^{\intts \fyr dr}\dzs d\widetilde{W}_s.
\end{align*}
Conditionning by $\F_t$ under $\widetilde{\mathbb{P}}$ and taking into account the fact that $Y_t$ and $\dyt$ are adapted with respect to $\F_t$ while $\intts \fyr dr$ and $\dzs$ are $\F_s$-adapted for $\theta < t \leq s \leq T$, we obtain 
\begin{align}
\label{EunEdeux}
\dyt =& \Ept{ e^{\intt \fys ds} D_{\theta}\xi }  \nonumber \\
&+ \Ept{\intt e^{\intts \fyr dr}\fxs\dxs ds}.
\end{align}
Using hypotheses $(\mbox{\textbf{H1}})$ and $(\mbox{\textbf{H2}})$, the first summand in (\ref{EunEdeux}) can be bounded by two positive constants $c$ and $C$ in the following manner:
\begin{align}
\label{premiersummand}
0 < c \leq \Ept{ e^{\intt \fys ds} D_{\theta}\xi } \leq C.
\end{align}
Using the results on $DX_{t}$ proved in Step 1 along with hypothesis $(\mbox{\textbf{H1}})$, we deduce that the second summand in (\ref{EunEdeux}) is bounded and non-negative: there is a positive constant $C$ such that
\begin{align}
\label{secondsummand}
0 \leq \Ept{\intt e^{\intts \fyr dr}\fxs\dxs ds} \leq C.
\end{align}
Combining the bounds (\ref{premiersummand}) and (\ref{secondsummand}), we immediately deduce that there exist two positive constants $c$ and $C$ such that
\begin{align}
\label{boundsurdy}
0 < c \leq D_{\theta}Y_{t} \leq C.
\end{align}
\\~\\
\textbf{\textit{Step 3: Conclusion of the proof by the Nourdin--Viens formula}}
\\~\\
Write $D_{\bullet}Y_t = \Phi_{Y_t}^{\bullet}(W)$ with a measurable function $\Phi_{Y_t}^{\bullet} : \mathbb{R}^{L^{2}([0,T])} \longrightarrow L^{2}([0,T])$. Then the bounds obtained in (\ref{boundsurdy}) yield, for $\theta < t$,
\begin{align*}
0 < c \leq \Phi_{Y_t}^{\theta}(W) \leq C.
\end{align*}
Define $\widetilde{\Phi_{Y_t}^{\bullet , u}}(W)= \Phi_{Y_t}^{\bullet}(e^{-u}W + \sqrt{1-e^{-2u}}W')$ for $u \in [0,+ \infty[$, where $W'$ stands for an independent copy of $W$ and is such that $W$ and $W'$ are defined on the product probability space $(\Omega \times \Omega', \mathcal{F}\otimes \mathcal{F}', \mathbb{P} \times \mathbb{P}')$. It is clear that, for $\theta < t$, we have
for any $u\in[0,\infty)$, $$0 < c \leq \widetilde{\Phi_{Y_t}^{\theta,u}}(W) \leq C.$$ Combining the two previous bounds yields, for $\theta < t$, $u\in[0,\infty)$, 
\begin{align}
\label{BornesProduitPhiYWPhiYWTilde}
0 < c^{2} \leq \Phi_{Y_t}^{\theta}(W)\widetilde{\Phi_{Y_t}^{\theta,u}}(W) \leq C^{2}.
\end{align}
Using the notation from Propositions \ref{existenceBornes} and \ref{PropPourCalcg}, we define
\begin{align*}
g(y)  =& \int_{0}^{\infty}e^{-u}\mathbf{E}\left(\mathbf{E'}\left( \int_{0}^{t}\Phi_{Y_t}^{\theta}(W)\widetilde{\Phi_{Y_t}^{\theta,u}}(W)d\theta \right) \Big|Y_t - \mathbf{E}(Y_t) = y\right)du.
\end{align*}
The bounds obtained in (\ref{BornesProduitPhiYWPhiYWTilde}) immediately yield
\begin{align*}
0 < ct \leq g(y) \leq Ct,
\end{align*}
with strictly positive constants $c$ and $C$. Thus, Propositions \ref{existenceBornes} and \ref{PropPourCalcg} conclude the proof of Theorem \ref{theoY}.
\qed \\

\section{Density of $Z_t$ : existence and Gaussian estimates}
We consider equations (\ref{eqdiff}) and (\ref{eqbackSDE}) with a function $f^{\star}$ that only has a linear dependency on $Z$, i.e.
\begin{numcases}
~X_t = x_0 + \into b(X_s) ds + \into \sigma(X_s) dW_s \label{eqdiffZ} \\
Y_t = \xi + \intt f^{\star} \xyzs ds - \intt Z_s dW_s      \label{eqbackSDEZ}
\end{numcases}
where $f^{\star}(x,y,z) = f(x,y) + \alpha z$, $\alpha \in \mathbb{R}$.
\\~\\
Because of the dependency of $f$ on $Z$, the Malliavin derivative $DZ$ will depend on $D^2Z$, which is not suitable for analyzing it within our framework. One can circumvent the above mentioned issue by using the Girsanov theorem to dispose of the impeding terms (similarly as done in the proof of Theorem \ref{theoY}). In order to clarify the proofs and to improve readability, we will consider that this step has already been performed in all of our proofs. This procedure leaves us with an equation of the type
\begin{numcases}
~X_t = x_0 + \into b(X_s) ds + \into \sigma(X_s) dW_s \label{eqdiffZ2var}\\
Y_t = \xi + \intt f \xys ds - \intt Z_s dW_s,       \label{eqbackSDEZ2var}
\end{numcases}
which is the one that will be referred to in the proofs of the upcoming results.
\\~\\ 
In the following subsections, we will prove our two main results concerning $Z_{t}$. We will begin by giving sufficient conditions for $Z_{t}$ to have a density. Up to our knowledge, this is the first result on density existence for the component $Z$ of the solution to equation (\ref{eqbackSDEZ}). We will then study in what framework and under what conditions this density can be bounded by Gaussian estimates.

\subsection{Existence of a density for $Z_{t}$}
We list in the next section the full set of hypotheses we need in this section.
\subsubsection{Hypotheses}\label{hyp-Z}
 We consider $b$, $\sigma$ and $f^{\star}$ to be appropriately smooth functions to ensure the existence and uniqueness of solutions to equations (\ref{eqdiffZ}) and (\ref{eqbackSDEZ}). We also impose additional conditions needed to prove Theorem \ref{theoZ}.
\begin{numcases}
~\mbox{\textbf{H4 : }} \xi \in L^{2}(\Omega, \mathcal{F}_{T})\cap \mathbb{D}^{2,2},\mbox{\ } \forall \theta \leq T,\mbox{\ } D_{\theta}\xi \geq 0 \mbox{\ } a.s  \mbox{\ and\ } \forall \mbox{\ } \theta < t \leq T,\mbox{\ } D_{\theta,t}^{2}\xi > 0 \mbox{\ } a.s.  \nonumber \\
\mbox{\textbf{H5 : }} f \in \mathcal{C}^{2}(\mathbb{R}^{2})   \mbox{\ and\ }  f_{x}, f_{y}, f_{xy}, f_{xx}, f_{yy}  \geq 0 \mbox{\ } a.s. \nonumber \\
\mbox{\textbf{H6 : }} b \in \mathcal{C}^{2}(\mathbb{R}), \mbox{\ } \sigma \in \mathcal{C}^{3}(\mathbb{R}), \mbox{\ } \sigma, \sigma', -\sigma'', -\sigma''' \geq 0 \mbox{\ } a.s \mbox{\ and\ }  [\sigma, [\sigma ,b]] \geq 0 \mbox{\ } a.s.  \nonumber
\end{numcases}
where $[b,\sigma]$ denotes the Lie bracket between $b$ and $\sigma$. 
\begin{remark}
Note that it is natural to have a condition on the iterated Lie bracket $[\sigma, [\sigma ,b]]$ between $b$ and $\sigma$ as second order Malliavin derivatives appear in the expression of $Z_{t}$.
\end{remark}
\begin{remark}
\label{remarksurlesigne}
The hypotheses on the signs of $\sigma$ and $f_{x}$ are made in order to make the proofs as readable as possible. It is possible to have more complex hypotheses for the signs of the products of derivatives of $\sigma$ and derivatives of $f$.
\end{remark}

\noindent The following theorem states that under the above hypotheses, $Z_{t}$ has a density on $\mathbb{R}$
\begin{theorem}
\label{theoZ}
Under the above hypotheses, for $t\in(0,T)$ the law of the random variable $Z_t$ is absolutely continuous with respect to the Lebesgue measure on $\mathbb{R}$.
\end{theorem}
\noindent Before proving Theorem \ref{theoZ}, we will first give a technical Lemma and a Proposition 
which will play a key role in the upcoming proof of this Theorem. 
First recall a lemma used to calculate the Malliavin derivative of a product of random variables in $\mathbb D^{1,2}$ 
(for example, see \cite{N}, p.36, exercice 1.2.12).
\begin{lemma}
\label{lem-nualart}(i)
Let $s, t\in[0,T]$ and $F\in\mathbb D^{1,2}$; then we have $\E\lp F|\F_t\rp\in\mathbb D^{1,2}$ and
\begin{align*}
D_s\E\lp F|\F_t\rp = \E\lp D_sF|\F_t\rp 1_{s\leq t}.
\end{align*}
(ii) If $F,G\in\mathbb D^{1,2}$ are such that $F$ and $\Vert DF\Vert _{L^{2}([0,T])}$ are bounded, then $FG\in \mathbb D^{1,2}$ and
\begin{align*}
D(FG)=FDG+GDF.
\end{align*}
\end{lemma}
\noindent The rest of this section is devoted to the proof of Theorem \ref{theoZ} which is divided in three steps. In the first step, we will prove that under the conditions of subsection \ref{hyp-Z} the second-order Malliavin derivatives of  $X$ and $Y$ are non-negative. This will be of importance in Step 2 (as these second-order derivatives appear in the expression of $DZ_{t}$.
\\~\\
In the second step, we will show that $DZ_{t}$ is positive almost surely. This will ensure that $\Vert DZ_{t}\Vert _{L^{2}([0,T])} > 0$ a.s.
\\~\\
In the third and last step, we will use the Bouleau--Hirsch Theorem to conclude the proof (see Theorem \ref{bouleauhirsch}).
\\~\\
\noindent {\bf Proof of Theorem \ref{theoZ}: }
\\~\\
\textbf{\textit{Step 1: Non-negativity of $D^{2}X$ and $D^{2}Y$}}
\\~\\
We start by proving that for $0 < \theta < t < s \leq T$, $\mathbb{P}$--a.s $D_{\theta,t}^{2} X_s$ is non-negative. 
\\~\\
Applying the Malliavin derivative to (\ref{relationEntreDXetDU2}) and using the second point in Lemma \ref{lem-nualart}, we deduce for $\theta, t \leq s \leq T$, since $U_s = g(X_s)$, 
\begin{align}
\label{exprD2Xs}
D_{\theta,t}^{2} X_s = & 
(\sigma\circ g^{-1})'(U_s)D_{\theta}U_s D_{t}U_s + (\sigma\circ g^{-1})(U_s)D_{\theta,t}^{2} U_s \nonumber \\
= & (\sigma' \sigma)(X_s)D_{\theta}U_s D_{t}U_s + \sigma(X_s)D_{\theta,t}^{2} U_s.
\end{align}
Hypothesis $(\mbox{\textbf{H6}})$ ensures that the term $(\sigma' \sigma)(X_s)D_{\theta}U_s D_{t}U_s$ is non-negative. It remains to prove that the second summand in (\ref{exprD2Xs}) is also non-negative. As $\sigma$ is non-negative, we focus on proving that $D_{\theta,t}^{2} U_s$ is too. Applying once again the Malliavin derivative operator to (\ref{derU}) and using the second point in Lemma \ref{lem-nualart}, we deduce for $\theta < t \leq s$,
\begin{align*}
D_{\theta,t}^{2} U_s =& \int_{t}^{s}(\beta \circ g^{-1})''(U_r)D_{t}U_{r}D_{\theta}U_{r}dr + \int_{t}^{s}(\beta \circ g^{-1})'(U_r)D_{\theta,t}^{2} U_rdr \\
=& \int_{t}^{s}e^{\int_{r}^{s}(\beta \circ g^{-1})'(U_v)dv}(\beta \circ g^{-1})''(U_r)D_{t}U_{r}D_{\theta}U_{r}dr \\
=& \int_{t}^{s}(\beta \circ g^{-1})''(U_r)D_{r}U_{s}D_{t}U_{r}D_{\theta}U_{r}dr.
\end{align*}
Further calculations yield the following expression
\begin{align*}
\pbg''(x) = & \left(\sigma \left(\frac{\left[\sigma, b \right]' }{\sigma} - \frac{\left[\sigma, b \right]\sigma'}{\sigma^{2}} \right) - \frac{1}{2}\left(\sigma '' \sigma \right)' \sigma   \right) \circ g^{-1}(x) \\
&= \lp\frac{[\sigma, [\sigma,b]]}{\sigma} -\undeux \lp\sigma''\sigma\rp'\sigma\rp\circ g^{-1}(x).
\end{align*}
Again, Hypothesis $(\mbox{\textbf{H6}})$ ensures that the term $\pbg''(x)$ is non-negative and thus that $D_{\theta,t}^{2} U_s$ is too. We have finally obtained that $D_{\theta,t}^{2} U_s \geq 0$ a.s.
\\~\\
We will now focus on $D_{\theta,t}^{2} Y_s$ and prove that for $0 < \theta < t < s \leq T$, $\mathbb{P}$--a.s it is also non-negative. Applying once more the Malliavin derivative operator to $D_{\theta}Y_s$ in (\ref{numeroDyt}) and using the second point in Lemma \ref{lem-nualart}, since $f$ does not depend on $Z$ we obtain, for $0 \leq \theta < t \leq s \leq T$,
\begin{align*}
\Dtht
Y_s =& \Dtht \xi -\int_s^T\Dtht Z_rdW_r \\
& + \int_s^T \Big\{ f_{xx}\xyr \Dthdeux X_r\Dtun X_r + f_x\xyr\Dtht X_r \\
& \qquad\qquad+ \fyxr\lp \Dthdeux Y_r\Dtun X_r  + \Dthdeux X_r\Dtun Y_r\rp \\
& \qquad\qquad +\fyyr\Dthdeux Y_r\Dtun Y_r + f_y\xyr\Dtht Y_r  \Big\}dr.
\end{align*}
Since $\Dtht Y_s$ solves a linear equation and is $\F_s$-measurable, we have that, for $0 \leq \theta < t \leq s \leq T$,
\begin{align}
\label{derYdeuxi}
\Dtht Y_s =& \E\Big( e^{\int_s^T f_{y}\lp X_r,Y_r\rp dr} \Dtht \xi|\F_s\Big) \nonumber \\
& + \E\Big(\int_s^T e^{\int_s^r f_{y}\lp X_u,Y_u\rp du} \Big\{ f_{xx}\xyr \Dthdeux X_r\Dtun X_r + f_x\xyr\Dtht X_r \nonumber\\
&\quad + \fyxr\lp\Dthdeux Y_r\Dtun X_r + \Dthdeux X_r\Dtun Y_r\rp 
+\fyyr\Dthdeux Y_r\Dtun Y_r \Big\}dr \Big|\F_s\Big).
\end{align}
Using hypotheses $(\mbox{\textbf{H4}})$ and $(\mbox{\textbf{H5}})$ along with the fact that for $0 < \theta < t < s \leq T$, $\mathbb{P}$--a.s $D_{\theta,t}^{2} X_s$ is non-negative, we obtain that for $0 < \theta < t < s \leq T$, $\mathbb{P}$--a.s $\Dtht Y_s$ is non-negative.
\\~\\
\textbf{\textit{Step 2: Positivity of $DZ_{t}$}}
\\~\\
Using a representation of $Z$, we compute its Malliavin derivative and show that under the hypotheses of Subsection \ref{hyp-Z}, it is almost surely positive. We begin by giving a representation of $Z$. We do not use the one from \cite{PP} in terms of gradient, that is $Z_t=\sigma\lp X_t\rp\lp \nabla X_t\rp^{-1}\nabla Y_t$, but rather use the fact that $Z_t$ can be represented by use of the Clark-Ocone formula. 
Indeed, by the uniqueness of the solution $(Y,Z)$, $Z_t$ can be written as
\begin{align}
\label{referencedeZt}
Z_t =\E \lp D_t \xi + D_t\int_0^Tf\xys ds \Big|\F_t\rp \in\mathbb D^{1,2}.
\end{align}
Using this fact, we get for $t\in[0,T]$
\begin{align*}
Z_t = & \E \lp D_t \xi + \int_t^T \left\lbrace f_x\xys D_tX_s + f_y\xys D_tY_s \right\rbrace ds |\F_t\rp.
\end{align*}
Let $\theta\leq t$. 
We use both points of Lemma \ref{lem-nualart} in order to calculate the first order Malliavin derivative of $Z_t$. This leads, for $\theta \leq t$:
\begin{align}
\label{deriveedezClarkOcone}
\Dt Z_t = & \E \Big( \Dtht \xi + \int_t^T \Big\{ \fxx \Dt X_sD_tX_s + \fyx \lp\Dt Y_sD_tX_s + \Dt X_sD_tY_s\rp \nonumber \\
& \qquad \quad + \fyy \Dt Y_sD_tY_s + f_x\xys\Dtt X_s + f_y\xys\Dtt Y_s \Big\} ds \Big| \F_t \Big).
\end{align}
Using Hypotheses $(\mbox{\textbf{H4}})$ and $(\mbox{\textbf{H5}})$ along with the results obtained in Step 1, we obtain that for $0 < \theta < t \leq T$, $\mathbb{P}-a.s$, $D_{\theta}Z_t > 0$.
\\~\\
\textbf{\textit{Step 3: Conclusion of the proof by the Bouleau--Hirsch Theorem}}
\\~\\
For all $t \leq T$, we have $$\Vert DZ_{t}\Vert _{L^{2}([0,T])}^{2} = \int_{0}^{T}(\Dt Z_t)^{2}d\theta .$$ Using the fact that for $0 < \theta < t \leq T$, $\mathbb{P}-a.s$, $D_{\theta}Z_t > 0$ proved in Step 2, we deduce that $\Vert DZ_{t}\Vert _{L^{2}([0,T])} > 0$ a.s. Applying the Bouleau--Hirsch Theorem (see Theorem \ref{bouleauhirsch}) concludes the proof.
\qed

\begin{remark}
Theorem \ref{theoZ} has been proven under a set of hypotheses (those of Subsection \ref{hyp-Z}) based on the fact that $\sigma$ is positive. The case where $\sigma$ is negative was included neither in the proof nor in the hypotheses for the sake of clarity and readability of the paper. However, as already mentioned in Remark \ref{remarksurlesigne}, this case can be addressed (without any further difficulties) by using the following transformations: $\sigma \rightarrow \tilde \sigma := - \sigma$ and $W \rightarrow \tilde W := -W$. After performing those tranformations, it suffices to consider $ (\tilde X, \tilde Y, \tilde Z) = (X, Y, -Z)$ to be the solution of
\begin{numcases}
~\tilde X_t = x_0 + \into b(\tilde X_s) ds + \into \tilde\sigma(\tilde X_s) d\tilde W_s \nonumber \\
\tilde Y_t = \xi + \int_t^T f\lp \tilde X_r, \tilde Y_r\rp dr -\int_t^T\tilde Z_rd\tilde W_r      \nonumber
\end{numcases}
This  brings the problem back to the set of hypotheses of Subsection \ref{hyp-Z} and it can be dealt with using the techniques presented in the last section.
\end{remark}

\subsection{Gaussian bounds for the density of $Z_{t}$}

In this section, we study a particular case of equations (\ref{eqdiffZ2var}) and (\ref{eqbackSDEZ2var}) and show that under proper conditions, the density of $Z_{t}$ can be bounded from above and below by Gaussian estimates. The backward equation we study is the following:
\begin{align}
Y_t = \phi \left(W_{T}\right)  + \intt f(Y_{s})ds - \intt Z_s dW_s,       \label{simplboundsback}
\end{align}
\subsubsection{Hypotheses}\label{hyp-bornesgaussZ}
 We consider $f$ to be an appropriately smooth function to ensure the existence and uniqueness of solutions to equation (\ref{simplboundsback}). We also impose additional conditions needed to prove Theorem \ref{theoZG}.
\begin{numcases}
~\mbox{\textbf{H7 : }} \phi \in \mathcal{C}_{b}^{2}(\mathbb{R}) \mbox{\ and\ } \phi'' \geq c > 0.  \nonumber \\
\mbox{\textbf{H8 : }} f \in \mathcal{C}_{b}^{2}(\mathbb{R}) \mbox{\ and\ } f', f'' \geq 0.  \nonumber
\end{numcases} 
The following theorem states that under the above hypotheses, $Z_{t}$ has a density that can be bounded from above and below by Gaussian estimates.
\begin{theorem}
\label{theoZG}
Under the above hypotheses, for $t\in(0,T)$ the random variable $Z_t$ defined in (\ref{simplboundsback}) has a density $\rho_{Z_t}$. Furthermore, there exist strictly positive constants $c$ and $C$ such that, for almost all $y\in\mathbb{R}$, $\rho_{Z_t}$ satisfies the following:
\begin{align*}
\frac{\mathbf{E}\vert Z_t -\mathbf{E}(Z_t)\vert}{2ct}\mathrm{exp}\left(-\frac{(z-\mathbf{E}(Z_t))^{2}}{2Ct}\right)  
\leq \rho_{Z_t}(z) \leq \frac{\mathbf{E}\vert Z_t -\mathbf{E}(Z_t)\vert}{2Ct}\mathrm{exp}\left(-\frac{(z-\mathbf{E}(Z_t))^{2}}{2ct}\right).
\end{align*}
\end{theorem}
\noindent {\bf Proof: } We will proceed in two steps, the first one being dedicated to proving that the Malliavin derivative of $Z_{t}$ is bounded and bigger than a positive constant. The second step will be to use the Nourdin--Viens formula to conclude the proof.
\\~\\
\textbf{\textit{Step 1: Boundedness and positivity of $DZ_{t}$}}
\\~\\
$Y_{t}$ being defined as in equation (\ref{simplboundsback}), equation (\ref{EunEdeux}) becomes 
\begin{align*}
\dyt = \Et{ e^{\intt f'(Y_{s}) ds} \phi'(W_{T})},
\end{align*}
and equation (\ref{derYdeuxi}) becomes
\begin{align*}
\Dtht Y_s =& \E\Big( e^{\int_s^T f'\lp Y_r\rp dr} \phi''(W_{T})|\F_s\Big) + \E\Big(\int_s^T e^{\int_s^r f'\lp Y_u\rp du} f''(Y_{r})\Dthdeux Y_r\Dtun Y_r dr|\F_s\Big).
\end{align*}
Using hypotheses $(\mbox{\textbf{H7}})$ and $(\mbox{\textbf{H8}})$, we obtain that $0 \leq \vert \dyt \vert \leq C$ and $0 < c \leq \Dtht Y_s  \leq C$. We finally compute $D_{\theta}Z_{t}$ from equation (\ref{deriveedezClarkOcone}) and we get
\begin{align*}
\Dt Z_t = & \E \Big( \phi''(W_{T}) + \int_t^T \Big\{f''(Y_{s}) \Dt Y_sD_tY_s + f'(Y_{s})\Dtt Y_s \Big\} ds \Big| \F_t \Big).
\end{align*}
Using Hypotheses $(\mbox{\textbf{H7}})$ and $(\mbox{\textbf{H8}})$ again, we finally get
\begin{align}
\label{gausssurZ}
0 < c \leq \Dt Z_t  \leq C.
\end{align}
~\\
\textbf{\textit{Step 2: Conclusion of the proof by the Nourdin--Viens formula}}
\\~\\
Write $D_{\bullet}Z_t = \Phi_{Z_t}^{\bullet}(W)$ with a measurable function $\Phi_{Z_t}^{\bullet} : \mathbb{R}^{L^{2}([0,T])} \longrightarrow L^{2}([0,T])$. Then (\ref{gausssurZ}) yields, for $\theta < t$, $0 < c \leq \Phi_{Z_t}^{\theta}(W) \leq C$. As previously done, define $\widetilde{\Phi_{Z_t}^{\bullet , u}}(W) = \Phi_{Z_t}^{\bullet}(e^{-u}W + \sqrt{1-e^{-2u}}W')$ for $u \in [0,+ \infty[$. Using (\ref{gausssurZ}), it is clear that, for $\theta < t$, we have for $u\in[0,+\infty)$, $0 < c \leq \widetilde{\Phi_{Z_t}^{\theta,u}}(W) \leq C$.
Combining the bounds on $\Phi_{Z_t}^{\theta}$ and $\widetilde{\Phi_{Z_t}^{\theta,u}}$ yields, for $\theta < t$ and $u\in[0,+\infty)$,
\begin{align}
\label{BornesProduitPhiZWPhiZWTilde}
0 < c^{2} \leq \Phi_{Z_t}^{\theta}(W)\widetilde{\Phi_{Z_t}^{\theta,u}}(W) \leq C^{2}.
\end{align}
Finally, let
\begin{align*}
g(z)  =& \int_{0}^{\infty}e^{-u}\mathbf{E}\left(\mathbf{E'}\big(\langle\Phi_{Z_t}^{\bullet}(W), \widetilde{\Phi_{Z_t}^{\bullet,u}}(W) \rangle_{L^{2}(\left[0,T\right])}\big)\Big|Z_t - \mathbf{E}(Z_t) = z\right)du \\
=& \int_{0}^{\infty}e^{-u}\mathbf{E}\left(\mathbf{E'}\left( \int_{0}^{t}\Phi_{Z_t}^{\theta}(W)\widetilde{\Phi_{Z_t}^{\theta,u}}(W)d\theta \right) \Big|Z_t - \mathbf{E}(Z_t) = z\right)du.
\end{align*}
The bounds obtained in (\ref{BornesProduitPhiZWPhiZWTilde}) immediatly yield $0 < ct \leq g(z) \leq Ct$. Thus, Proposition \ref{existenceBornes} concludes the proof of Theorem \ref{theoZG}.
\qed 
\begin{remark}
It is also possible to derive Gaussian density estimates for more complex equations than the one dealt with in this section. Hypotheses have to be changed in each case, making it difficult to state a general result with reasonable hypotheses covering most cases.
\end{remark}
\subsection*{Acknowledgments:}
We would like to thank F.G. Viens for introducing us to this topic
as well as A. Millet and C.A. Tudor for helpful comments. \\~\\We also thank an anonymous referee for his/her thorough review and highly appreciate the comments and
suggestions, which significantly contributed to improving the quality of this paper.

\end{document}